\documentclass{article}
\usepackage[latin1]{inputenc}

\usepackage{abstract}
\usepackage{amsthm}
\newtheorem{theorem}{Theorem}

\usepackage{amssymb,amsmath}
\begin{document}
\title{Non-existence of CR submanifolds of maximal CR dimension satisfying $RA=0$ in non-flat complex space forms }
\author{Mirjana Milijevi\'{c}\\\\
Faculty of Architecture and Civil Engineering\\
University of Banja Luka\\
Stepe Stepanovi\'{c}a 77, 78000 Banja Luka\\
Bosnia and Herzegovina\\
E-mail: \texttt{$mirjana_{-}milijevic@yahoo.com$}}
\date{}
\maketitle
\begin{abstract}
It has been proved that there are no real hypersurfaces satisfying $RA=0$ in non-flat complex space forms. In this paper we prove that the same is true in the case of CR submanifolds of maximal CR dimension, that is there are no CR submanifolds of maximal CR dimension satisfying $RA=0$ in non-flat complex space forms.
\end{abstract}
\emph{Key words and phrases.} Complex space form, $CR$ submanifold of maximal $CR$ dimension, shape operator, curvature tensor.\\
\emph{AMS Subject Classification.} 53C15, 53C40, 53B20.
\section{Introduction}
S. Maeda proved in [4] the non-existence of real hypersurfaces satisfying $RA=0$ in the complex projective space, where we denoted by $R$ the curvature tensor and by $A$ the shape operator of a hypersurface. On the other hand M. Ortega proved in [2] that there are no real hypersurfaces in non-flat complex space forms such that $RA=0$.

As a real hypersurface is a typical example of a CR submanifold of maximal CR dimension, we will in this paper generalize the results obtained by S. Maeda and M. Ortega to CR submanifolds of maximal CR dimension.

Let $\overline{M}$ be an $(n + p)$-dimensional complex space form, i.e. a Kaehler manifold of constant holomorphic sectional curvature $4c$, endowed with metric $\overline{g}$. Let $M$ be an $n$-dimensional real submanifold of $\overline{M}$ and $J$ be the almost complex structure of $\overline{M}$. For a tangent space $T_{x}(M)$ of $M$ at $x$,\:we put $H_{x}(M)=JT_{x}(M)\cap T_{x}(M)$. Then,\:$H_{x}(M)$ is the maximal complex subspace of $T_{x}(M)$ and is called the holomorphic tangent space to $M$ at $x$. If the complex dimension $dim_{\textbf{C}}H_{x}(M)$ is constant over $M$,\:$M$ is called a Cauchy-Riemann submanifold or briefly a CR submanifold and the constant $dim_{\textbf{C}}H_{x}(M)$ is called the CR dimension of $M$. If,\:for any $x\in M$,\:$H_{x}(M)$ satisfies $dim_{\textbf{C}}H_{x}(M)=\frac{n-1}{2}$,\:then $M$ is called a CR submanifold of maximal CR dimension. It follows that there exists a unit vector field $\xi$ normal to $M$ such that $JT_{x}(M)\subset T_{x}(M)\oplus span\{\xi_{x}\}$,\:for any $x\in M$.

\section{$CR$ submanifolds of maximal $CR$ dimension of a complex space form}
Let $\overline{M}$ be an $(n + p)$-dimensional complex space form with Kaehler structure $(J,\overline{g})$ and of constant holomorphic sectional curvature $4c$. Let $M$ be an $n$-dimensional $CR$ submanifold of maximal $CR$ dimension in $\overline{M}$ and
$\iota :M\to {\overline{M}}$  immersion. Also,\:we denote by $\iota$ the differential of the immersion. The Riemannian metric $g$ of $M$ is induced from the Riemannian metric $\overline{g}$ of $\overline{M}$ in such a way that
$g(X,Y)=\overline{g}(\iota X,\iota Y)$,\:where $X,\:Y \in T(M)$. We denote by $T(M)$ and $T^{\bot}(M)$ the tangent bundle and the normal bundle of $M$,\:respectively.

On $\overline{M}$ we have the following decomposition into tangential and normal components:
\begin{align}\label{eq:m1}
J\iota X=\iota FX+u(X)\xi,\:\: X\in T(M).
\end{align}
Here $F$ is a skew-symmetric endomorphism acting on $T(M)$ and $u$ in one-form on $M$.

Since $T_{1}^{\bot}(M)=\{\eta\in T^{\bot}(M)|\overline{g}(\eta,\xi)=0\}$ is $J$-invariant,\:from now on we will denote the orthonormal basis of $T^{\bot}(M)$ by $\xi,\xi_{1},\cdots,\xi_{q},\xi_{1^{*}},\cdots,\xi_{q^{*}}$,\:where $\xi_{a^{*}}=J\xi_{a}$ and $q=\frac{p-1}{2}$. Also,\:$J\xi$ is the vector field tangent to $M$ and we write

\begin{align}\label{eq:m2}
J\xi=-\iota U.
\end{align}
Furthermore,\:using (\ref{eq:m1}),\:(\ref{eq:m2}) and the Hermitian property of $J$ implies
\begin{align}\label{eq:80}
F^{2}X=-X+u(X)U,
\end{align}
\begin{align}\label{eq:m5}
FU=0,
\end{align}
\begin{align}\label{eq:m12}
\nabla_{X}U = FAX,
\end{align}
\begin{align}\label{eq:m6}
g(X,U)=u(X).
\end{align}
Next,\:we denote by $\overline{\nabla}$ and $\nabla$ the Riemannian connection of $\overline{M}$ and $M$,\:respective-\-ly, and by $D$ the normal connection induced from $\overline{\nabla}$ in the normal bundle of $M$. They are related by the following Gauss equation
\begin{align}\label{eq:m7}
\overline{\nabla}_{\iota X}\iota Y=\iota \nabla_{X}Y+h(X,Y),
\end{align}
where $h$ denotes the second fundamental form,\:and by Weingarten equations
\begin{align}\label{eq:m8}
\overline{\nabla}_{\iota X}\xi&=-\iota AX+D_{X}\xi\\
\notag
&=-\iota AX+\sum_{a=1}^{q}\{s_{a}(X)\xi_{a}+s_{a^{*}}(X)\xi_{a^{*}}\},
\end{align}
\begin{align}\label{eq:m9}
\overline{\nabla}_{\iota X}\xi_{a}&=-\iota A_{a}X+D_{X}\xi_{a}=-\iota A_{a}X-s_{a}(X)\xi\\
\notag
&+\sum_{b=1}^{q}\{s_{ab}(X)\xi_{b}+s_{ab^{*}}(X)\xi_{b^{*}}\},
\end{align}
\begin{align}\label{eq:m10}
\overline{\nabla}_{\iota X}\xi_{a^{*}}&=-\iota A_{a^{*}}X+D_{X}\xi_{a^{*}}=-\iota A_{a^{*}}X-s_{a^{*}}(X)\xi\\
\notag
&+\sum_{b=1}^{q}\{s_{a^{*}b}(X)\xi_{b}+s_{a^{*}b^{*}}(X)\xi_{b^{*}}\},
\end{align}
where the $s$'s are the coefficients of the normal connection $D$ and $A$,\:$A_{a}$,\:$A_{a^{*}}$;\:$a=1,\cdots,q$,\:are the shape operators corresponding to the normals $\xi$,\:$\xi_{a}$,\:$\xi_{a^{*}}$,\:respecti-\-vely. They are related to the second fundamental form by
\begin{align}\label{eq:m11}
h(X,Y)&=g(AX,Y)\xi\\
\notag
&+\sum_{a=1}^{q}\{g(A_{a}X,Y)\xi_{a}+g(A_{a^{*}}X,Y)\xi_{a^{*}}\}.
\end{align}
Since the ambient manifold is a Kaehler manifold,\:using (\ref{eq:m1}),\:(\ref{eq:m2}),\:(\ref{eq:m9}) and (\ref{eq:m10}),\:it follows that
\begin{align}\label{eq:m14}
s_{a^{*}}(X)=u(A_{a}X),
\end{align}
\begin{align}\label{eq:m15}
s_{a}(X)=-u(A_{a^{*}}X),
\end{align}
for all $X,\:Y$ tangent to $M$ and $a=1,\cdots,q$.\\
The Codazzi and the Gauss equation for the distinguished vector field $\xi$ are
\begin{align}\label{eq:m21}
&(\nabla_{X}A)Y-(\nabla_{Y}A)X=c\{u(X)FY-u(Y)FX-2g(FX,Y)U\}\\
\notag
&+\sum_{a=1}^{q}\{s_{a}(X)A_{a}Y-s_{a}(Y)A_{a}X\}+\sum_{a=1}^{q}\{s_{a^{*}}(X)A_{a^{*}}Y-s_{a^{*}}(Y)A_{a^{*}}X\},
\end{align}
\begin{align}\label{eq:m22}
R_{XY}Z&=c\{g(Y,Z)X-g(X,Z)Y+g(FY,Z)FX\\
\notag
&-g(FX,Z)FY - 2g(FX,Y)FZ\}\\
\notag
&+g(AY,Z)AX-g(AX,Z)AY\\
\notag
&+\sum_{a=1}^{q}\{g(A_{a}Y,Z)A_{a}X-g(A_{a}X,Z)A_{a}Y\}\\
\notag
&+\sum_{a=1}^{q}\{g(A_{a^{*}}Y,Z)A_{a^{*}}X-g(A_{a^{*}}X,Z)A_{a^{*}}Y\},
\end{align}
respectively, for all $X$, $Y$, $Z$ tangent to $M$, where $R$ denotes the Riemannian curvature tensor of $M$.

\section{CR submanifolds of maximal CR dimension satisfying $RA=0$}
\begin{theorem}\label{T1}
Let $M$ be an $n$-dimensional CR submanifold of maximal CR dimension in an $(n+p)$-dimensional complex space form
$(\overline{M},J,\overline{g})$, where $n\geq 3$ and the constant holomorphic sectional curvature of $\overline{M}$ equals $4c$.
Let $p < n$, $A$ be the shape operator of the distinguished vector field $\xi$ and $R$ be the Riemannian curvature tensor of $M$. If $RA=0$ on $M$, then $\overline{M}$ is an Euclidean space.
\end{theorem}
\flushleft \emph{Proof.}
Because of the assumption that $RA=0$ we have
\begin{align}
\notag
g(R_{XY}(AZ),W) = g(AR_{XY}Z,W),
\end{align}
for $X$, $Y$, $Z$, $W$ tangent to $M$, that is
\begin{align}\label{eq:m23}
c\{&g(Y,AZ)g(X,W)-g(X,AZ)g(Y,W)+g(FY,AZ)g(FX,W)\\
\notag
&-g(FX,AZ)g(FY,W)-2g(FX,Y)g(FAZ,W)\}\\
\notag
& + g(AY,AZ)g(AX,W) - g(AX,AZ)g(AY,W)\\
\notag
&+\sum_{a=1}^{q}\{g(A_{a}Y,AZ)g(A_{a}X,W) - g(A_{a}X,AZ)g(A_{a}Y,W)\}\\
\notag
&+\sum_{a=1}^{q}\{g(A_{a^{*}}Y,AZ)g(A_{a^{*}}X,W) - g(A_{a^{*}}X,AZ)g(A_{a^{*}}Y,W)\}=\\
\notag
c\{&g(Y,Z)g(AX,W) - g(X,Z)g(AY,W) + g(FY,Z)g(AFX,W)\\
\notag
& - g(FX,Z)g(AFY,W)-2g(FX,Y)g(AFZ,W)\} \\
\notag
&+ g(AY,Z)g(A^{2}X,W) - g(AX,Z)g(A^{2}Y,W)\\
\notag
&+\sum_{a=1}^{q}\{g(A_{a}Y,Z)g(AA_{a}X,W) - g(A_{a}X,Z)g(AA_{a}Y,W)\}\\
\notag
&+\sum_{a=1}^{q}\{g(A_{a^{*}}Y,Z)g(AA_{a^{*}}X,W) - g(A_{a^{*}}X,Z)g(AA_{a^{*}}Y,W)\}
\end{align}
Interchanging $X$ and $Z$ in (\ref{eq:m23}) and subtracting the resulting equation and (\ref{eq:m23}) we obtain
\begin{align}\label{eq:m24}
c\{&-g(FX,AZ)g(FX,Z) - g(X,AX)g(Z,Z) + 3g(FZ,X)g(FAX,Z)\\
\notag
&+g(X,X)g(Z,AZ) - 4g(FZ,X)g(FX,AZ)\}\\
\notag
&-g(AX,AX)g(AZ,Z) + g(AX,X)g(AZ,AZ)\\
\notag
&+\sum_{a=1}^{q}\{g(A_{a}X,X)g(A_{a}Z,AZ) - g(A_{a}X,AZ)g(A_{a}X,Z)\}\\
\notag
&+\sum_{a=1}^{q}\{g(A_{a^{*}}X,X)g(A_{a^{*}}Z,AZ) - g(A_{a^{*}}X,AZ)g(A_{a^{*}}X,Z)\}\\
\notag
&+\sum_{a=1}^{q}\{g(A_{a}Z,AX)g(A_{a}X,Z) - g(A_{a}X,AX)g(A_{a}Z,Z)\}\\
\notag
&+\sum_{a=1}^{q}\{g(A_{a^{*}}Z,AX)g(A_{a^{*}}X,Z) - g(A_{a^{*}}X,AX)g(A_{a^{*}}Z,Z)\}=0.
\end{align}
From (\ref{eq:m24}) it follows
\begin{align}\label{eq:m25}
c\{&-g(X,AX)Z  3g(FAX,Z)FX + g(X,X)AZ + 3g(FX,AZ)FX\}\\
\notag
&-g(AX,AX)AZ + g(AX,X)A^{2}Z\\
\notag
+\sum_{a=1}^{q}\{&g(A_{a}X,X)AA_{a}Z - g(A_{a}X,Z)AA_{a}X + g(A_{a^{*}}X,X)AA_{a^{*}}Z\\
\notag
&-g(A_{a^{*}}X,Z)AA_{a^{*}}X + g(A_{a}X,Z)A_{a}AX - g(A_{a}X,AX)A_{a}Z\\
\notag
&+g(A_{a^{*}}X,Z)A_{a^{*}}AX - g(A_{a^{*}}X,AX)A_{a^{*}}Z\}=0,
\end{align}
because $Z$ is an arbitrary tangent vector.\\
On the other hand, from (\ref{eq:m24}) it follows
\begin{align}\label{eq:m26}
c\{&-g(X,AX)Z + 3g(FZ,X)FAX + g(X,X)AZ - 3g(FZ,X)AFX\}\\
\notag
&-g(AX,AX)AZ + g(AX,X)A^{2}Z\\
\notag
+\sum_{a=1}^{q}\{&g(A_{a}X,X)AA_{a}Z- g(A_{a}X,Z)AA_{a}X + g(A_{a^{*}}X,X)AA_{a^{*}}Z\\
\notag
&-g(A_{a^{*}}X,Z)AA_{a^{*}}X + g(A_{a}X,Z)A_{a}AX - g(A_{a}X,AX)A_{a}Z\\
\notag
&+g(A_{a^{*}}X,Z)A_{a^{*}}AX - g(A_{a^{*}}X,AX)A_{a^{*}}Z\}=0,
\end{align}
because $Z$ is an arbitrary tangent vector.\\
Subtracting (\ref{eq:m25}) and (\ref{eq:m26}) we obtain
\begin{align}\label{eq:m27}
c\{-g(FAX,Z)FX + g(FX,AZ)FX - g(FZ,X)FAX + g(FZ,X)AFX\}=0.
\end{align}
After putting $Z=X$ in (\ref{eq:m27}) we obtain
\begin{align}\label{eq:m28}
c\{-g(FAX,X)FX + g(FX,AX)FX\}=0.
\end{align}
Multiplying (\ref{eq:m28}) with $AX$, we obtain
\begin{align}\label{eq:m29}
2cg(FAX,X)^{2} = 0.
\end{align}
If $c\neq 0$, then from (\ref{eq:m29}) it follows that
\begin{align}\label{eq:m30}
FAX=0.
\end{align}
From (\ref{eq:m30}) we conclude that
\begin{align}\label{eq:m31}
AX=\alpha U,
\end{align}
for some function $\alpha$ and $X\in T(M)$.\\
Multiplying the Codazzi (\ref{eq:m21}) equation with $U$ and putting $Y=U$, we obtain
\begin{align}\label{eq:m32}
g((\nabla_{X}A)U,U) &= g((\nabla_{U}A)X,U) \\
\notag
&+ \sum_{a=1}^{q}\{2g(A_{a^{*}}U,U)g(A_{a}U,X)-2g(A_{a^{*}}U,X)g(A_{a}U,U)\}.
\end{align}
Differentiating (\ref{eq:m31}), we obtain
\begin{align}\label{eq:m33}
(\nabla_{Y}A)X = (Y\alpha)U - \alpha U,
\end{align}
$Y\in T(M)$.\\
Now, from (\ref{eq:m32}) and (\ref{eq:m33}) we get
\begin{align}\label{eq:m34}
X\alpha = U\alpha + \sum_{a=1}^{q}\{2g(A_{a^{*}}U,U)g(A_{a}U,X)-2g(A_{a^{*}}U,X)g(A_{a}U,U)\}.
\end{align}
From the Codazzi (\ref{eq:m21}) equation multiplied with $U$, (\ref{eq:m33}) and (\ref{eq:m34}) it follows
\begin{align}\label{eq:m35}
&\sum_{a=1}^{q}\{g(A_{a^{*}}U,U)g(A_{a}U,X) - g(A_{a}U,U)g(A_{a^{*}}U,X)\} +\\
\notag
&\sum_{a=1}^{q}\{-g(A_{a^{*}}U,U)g(A_{a}U,Y) + g(A_{a}U,U)g(A_{a^{*}}U,Y)\}=\\
\notag
&-cg(FX,Y) + \sum_{a=1}^{q}\{g(A_{a^{*}}U,Y)g(A_{a}U,X) - g(A_{a^{*}}U,X)g(A_{a}U,Y)\}
\end{align}
Putting $X=FX$ in (\ref{eq:m35}) and then $X=U$ in the resulting equation, we obtain
\begin{align}\label{eq:m36}
\sum_{a=1}^{q}\{-g(A_{a^{*}}U,U)g(A_{a}U,Y) + g(A_{a}U,U)g(A_{a^{*}}U,Y)\}=0.
\end{align}
From (\ref{eq:m35}) and (\ref{eq:m36}) we obtain
\begin{align}\label{eq:m37}
cFX = \sum_{a=1}^{q}\{g(A_{a}U,X)A_{a^{*}}U - g(A_{a^{*}}U,X)A_{a}U\}.
\end{align}
From (\ref{eq:m37}) it follows that $FX$ is a linear combination of $A_{a}U$ and $A_{a^{*}}U$; $a=1,\cdots,q$.\\
Since every tangent vector $Y$ orthogonal to $U$ can be expressed as $Y=FX$, for $n-1>2q=p-1$, i.e. $n>p$, it follows that there exists a unit vector field $Y=FX$ which is orthogonal to $span\{A_{a}U,A_{a^{*}}U\}$; $a=1,\cdots,q$.\\
Putting such $Y=FX$ in (\ref{eq:m37}) instead of $X$ we get
\begin{align}\label{eq:m38}
cF^{2}X=0,
\end{align}
from which it follows $F^{2}X=0$, for $X\in T(M)$.\\
This is a contradiction because of (\ref{eq:80}). $\Box$
\begin{theorem}\label{T2}
Let $\textbf{M}$ be an $n$-dimensional $CR$ submanifold of maximal $CR$ dimension in an $(n+p)$-dimensional complex space form $(\overline{\textbf{M}},J,\overline{g})$, where $n\geq 3$ and the constant holomorphic sectional curvature of $\overline{\textbf{M}}$ equals $4c$.  Let the distinguished vector field $\xi$ be parallel with respect to the normal connection $D$, $A$ be the shape operator of $\xi$ and $R$ be the Riemannian curvature tensor of $M$. If $RA=0$ on $M$, then $\overline{M}$ is an Euclidean space.
\end{theorem}
\flushleft \emph{Proof.}
As in the proof of the Theorem 1 we obtain
\begin{align}\label{eq:m39}
(\nabla_{Y}A)X=(Y\alpha)U - \alpha U,
\end{align}
for some function $\alpha$ and $X$, $Y\in T(M)$.\\
Putting $Y=U$ in the Codazzi equation (\ref{eq:m21}) and multiplying the resulting equation with $U$ we obtain
\begin{align}\label{eq:m40}
g((\nabla_{X}A)U - (\nabla_{U}A)X,U)=0.
\end{align}
From (\ref{eq:m39}) and (\ref{eq:m40}) it follows
\begin{align}\label{eq:m41}
X\alpha=U\alpha.
\end{align}
From the the Codazzi equation (\ref{eq:m21}) multiplied with $U$ and (\ref{eq:m41}) we obtain
\begin{align}\label{eq:m42}
-2cg(FX,Y)=0,
\end{align}
from which it follows that $g(FX,Y)=0$.\\
This is a contradiction. $\Box$

\end{document}